\documentclass[]{amsart}

\newtheorem{theorem}{Theorem}[section]

\newtheorem{lemma}[theorem]{Lemma}
\newtheorem{corollary}[theorem]{Corollary}
\theoremstyle{definition}
\newtheorem{definition}[theorem]{Definition}
\theoremstyle{remark}

\numberwithin{equation}{section}

\usepackage{amssymb}
\usepackage{graphicx}
\usepackage[citecolor=blue,linkcolor=blue,linkbordercolor={1 1 1},colorlinks=true,citebordercolor={1 1 1}]{hyperref}
\usepackage{pstricks}
\usepackage{pst-node}

\begin{document}

\title{Graph-Chromatic Implicit Relations}

\author{JOSE ANTONIO MARTIN H.}
\address{Dept. Sistemas Inform\'aticos \&
Computaci\'on, Universidad Complutense de Madrid (Spain). }
\curraddr{Dept. Sistemas Inform\'aticos \& Computaci\'on,
Universidad Complutense de Madrid, Ciudad Universitaria, Fac. de
Inform\'atica, C. Prof. Jos\'e Garc\'ia Santesmases, s/n. 28040
(Spain) } \email{jamartinh@fdi.ucm.es}





\date{Oct. 2, 2006 and, in revised form, -, -}



\begin{abstract}
A theory about the implication structure in graph coloring is
presented. Discovering hidden relations is a crucial activity in
every scientific discipline. The development of mathematical
models to study and discover such hidden relations is of the most
highest interest. The main contribution presented in this work is
a model of hidden relations materialized as implicit-edges and
implicit-identities in the graph coloring problem, these relations
can be interpreted in physical and chemical models as hidden
forces, hidden interactions, hidden reactions or hidden variables.
Also this theory can be extended to the complete class of
NP-complete problems.

\end{abstract}

\maketitle
\tableofcontents

\section{Introduction}

All the graphs in this work are simple, that is, unweighted, undirected and containing no graph loops or multiple edges.

Graph Vertex Coloring is the task to assign a color to every vertex of a graph in such a way that no two adjacent vertices are
assigned the same color, and the chromatic number of a graph, is the minimum number of distinct colors needed to properly color a
graph. The importance of graph coloring resides in its general applicability as a process of partitioning a collection of objects
into disjoint sets fulfilling a set of constraints between that objects. For a good reference on graph coloring problems the
reader can review the excellent book of Jensen and Toft~\cite{JT95}. We are presenting a theory about the implication
structure in graph coloring guided by the concept of implicit edge. We think that the study of this kind of relations has a high
applicability in a wide range of experimental and theoretical fields of research. This study can help to unveil a wide range of
``hidden" phenomena in a wide range of scientific disciplines like physics, chemistry, and of course, discrete mathematics and
computer science.

In this work, a formal definition of the Theory of graph-chromatic implicit-relations is presented. The notion of implicit-relation is taken in this work as any relation between a set of vertices of a graph such that despite those relations are not contained explicitly in any set of constraints or relations set, for instance the set of edges of a graph, such relations still exist for the graph vertex coloring problem. In this way we will define here the implicit-edge relation which means that there is a kind of invisible (but real) edge between two vertices and implicit-identities which means that two different vertices could be treated as the same vertex.

In order to formalize such concepts a series of new lemmas and theorems are presented which relate directly these novel concepts with some key aspects of the current theory of graph coloring such as critical graphs, graphs minors, Kempe chains and also with interesting open problems such as the Hadwiger conjecture and the double-critical graphs problem.

\section{Preliminary Definitions}
Partitioning the set of vertices $V(G)$ of a graph $G$ into separate classes, in such a way that no two adjacent vertices are
grouped into the same class, is called the vertex graph coloring problem, In order to distinguish such classes, a set of colors C
is used, and the division into these classes is given by a `proper' coloring $\varphi : V(G)\rightarrow C$, where
$\varphi(u)\neq \varphi(v)$ for all $(u,v)$ belonging to the set of edges $E(G)$ of $G$. An independent set $S$ of a graph $G$ is a subset of vertices of $G$ such that there is no edge between any pair of vertices of $S$ in $E(G)$. Hence we can also define the graph coloring problem as the task of separating the vertices of a graph into a collection of different independent sets.

Given a graph vertex coloring problem over a graph $(G)$ with a set of colors $(C)$, if  $(C)$  has cardinality  $(k)$ , then
$\varphi$ is a \emph{k-coloring} of  $(G)$. The \emph{Chromatic number} of a graph $\chi(G)$ is the minimum number of colors
necessary to color the vertices of a graph $G$ in such a way that no two adjacent vertices are colored with the same color, thus, if
$\chi(G)\leq k$ then one says that G is \emph{k-colorable} and if $\chi(G)= k$ then one says that G is \emph{k-chromatic}. The
neighborhood $N(v)$ of $v\in V(G)$ is the complete set of vertices that are connected to $v$. An interesting problem in graph coloring is the
counting task \#k-coloring, this problem consists in counting the total number $P(G,k)$ of distinct k-colorings of a given graph G, it is believed that in computational complexity terms this counting problem is harder than the decision problem of just determining k-colorability in the general case. One standard way of solving the \#k-coloring is by using the so called chromatic polynomial which is a polynomial that encodes the number $P(G,k)$.


A minor of a graph $G$ is any graph obtained by successive contraction or deletion of edges.

Let $\mathcal{S^*}(G)=\{S_1,...,S_n\}$ be the superset of all independent sets of a graph $G$. Let us consider the next four
operations that could affect $\mathcal{S^*}(G)$:
\begin{enumerate}
\item to draw a new edge $\{x,y\}$ in $G$: $(G + xy)$.
\item to delete a vertex $x$ of $G$: $(G-x)$.
\item to identity a pair of vertices $\{x,y\}$ of $G$: $(G/xy)$.
\item to contract an edge $e=\{x,y\}$ of $G$: $(G/e)$.
\end{enumerate}

The details of each operation is shown next:
\begin{enumerate}
\item  It is easy to see that the operation of drawing a new edge could not produce any new independent set, indeed it
eliminates at least one independent set from $\mathcal{S^*}$.
\item  It is also easy to see that the operation of deleting a vertex $x$ from $G$ could not produce any new
independent set, indeed it eliminates at least as many sets from $\mathcal{S^*}$ as the number of sets in which $x$ appears.
\item and 4. Now we show that the two operations (identify and contract) does not produce any new independent set by
decomposing them in the two operations firstly studied.
\end{enumerate}
Both, vertex identification and edge contraction could be realized by the next simple procedure: 
\begin{itemize}
\item Draw a new edge from any vertex joined to $y$ towards $x$ and then delete vertex $y$.
\end{itemize}
As we can see, both operations drawing a new edge and delete a vertex could not produce any new independent
set. This fact can be formalized as the next lemma:

\begin{lemma}
\label{IS} Given a graph G and pair of vertices $\{x,y\}$ of $G$.

Let $\mathcal{S^*}(G)=\{S_1,...,S_n\}$ be set superset of all
independent sets of $G$ and let $H$ be the graph obtained by
contracting and edge of $G$ or by identifying a pair of vertices
$\{x,y\}$ of $G$ and respectively let
$\mathcal{S^*}(H)=\{S_1,...,S_n\}$ be the superset of all
independent sets of $H$ then:
\begin{equation}
\mathcal{S^*}(H)\subset \mathcal{S^*}(G)
\end{equation}
\end{lemma}

Hence any of the four operations listed does not increment the number of independent sets of a graph.

\section{Implicit Relations} \label{sec:12} The chromatic
implicit relations are mainly defined by their two most basics
concepts, implicit-edges and implicit-identities.

\subsection{Implicit Edges}

The notion of implicit-edge is based on the simple idea that in
the problem of graph coloring with a fixed number of colors (k),
that is, to find a k-coloring of some k-chromatic graph G, there
are some independent sets $S$ of vertices of $G$ such that despite those
vertices are not joined by an edge we can not find any k-coloring
of G where all the vertices in $S$ are colored with the same color. In order to
define precisely the concept of implicit-edge we also fix the
cardinality of the subsets $S$ to $|S|=2$, that is, we define an
implicit-edge as a 2-subset $\{i,j\}$ of vertices of G.

Also, it is shown that the notion of implicit-edge is invariant to
the fact that the 2-subset $\{i,j\}$ belongs to the set of edges
of G or not. If a 2-subset $\{i,j\}$ is an implicit-edge and
$\{i,j\}$ belongs to the set of edges of G we say that $\{i,j\}$
is a drawn~implicit-edge and if $\{i,j\}$ does not belongs to the
set of edges of G we say that $\{i,j\}$ is a non-drawn~implicit-edge.

Now we will define in strict sense the notion of implicit-edge:

\begin{definition}Implicit Edge.\\
\label{IE} Given a graph G and an integer k, we say that the
subset $\{i,j\}$ of vertices of G is an implicit-edge \textbf{if
and only if there is no} k-coloring of the graph $G-ij$ such that
vertices $i$ and $j$ receive the same color being the graph $G-ij$
the graph G with a possible edge $\{i,j\}$ deleted, hence
$\{i,j\}$ is an implicit-edge in the graph $G-ij$ too.
\end{definition}

Alteratively, we can say that the set of all k-colorings of $G-ij$
where $i$ and $j$ receive the same color is the empty set.

Also we must note that implicit-edges $\{i,j\}$ could belong to the set
of edges of G or not but the set of all k-colorings is the same in both cases.

In order to illustrate this notion we will use a trivial case,
that is, a bipartite (2-chromatic) graph G. Implicit-edges in
2-chromatic graphs are defined very easy: every subset $\{i,j\}$
of vertices of G with an odd length path between them are
implicit-edges. Let us take the graph G as the four-path:
$$\mbox{$i\circ$---$\circ$---$\circ$---$\circ j$}$$

We can easily see that $\{i,j\}$ is an implicit-edge of G indeed
the only implicit-edge of G. Also, if we draw the edge $\{i,j\}$
we get a square and then all edges are now implicit-edges.

More complex examples are shown next for the 3-coloring case. In the graph (a) of Fig.~\ref{impl1}., there is an implicit edge
between the vertices \{$u$,$v$\}. As one can easily see, coloring the
vertices \{$u$,$v$\} with the same color leads to an impossibility for
coloring the graph with only three colors. The same fact can be
verified in the graph (b) of Fig.~\ref{impl1}., there is an
implicit edge between the vertices \{$u$,$v$\}.

\newlength{\MyLength}
\settowidth{\MyLength}{$2$}
\newcommand{\MyBox}[1]{\makebox[\MyLength]{\tiny#1}}

\begin{figure}[tbh]
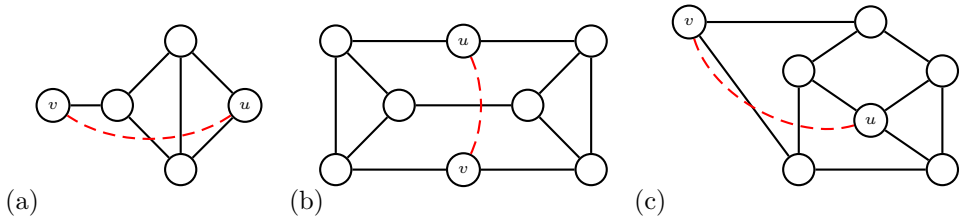

\begin{centering}
  \begin{tabular}{lll}
   \begin{psmatrix}[colsep=0.4,rowsep=0.4,mnode=circle]
  &&& [name=A1]\MyBox{$$} \\
  &[name=A5]\MyBox{$v$}& [name=A3]\MyBox{$$} && [name=A2]\MyBox{$u$} \\
  &&& [name=A4]\MyBox{$$}
  \psset{nodesep=0pt,arrows=-}
  \ncline{A1}{A2}
  \ncline{A1}{A3}
  \ncline{A1}{A4}
  \ncline{A2}{A4}
  \ncline{A3}{A4}
  \ncline{A3}{A5}
  \ncarc[arcangle=35,linecolor=red,linestyle=dashed]{A2}{A5}
 \end{psmatrix}
 &
 \begin{psmatrix}[colsep=0.4,rowsep=0.4,mnode=circle]
  &[name=A1]\MyBox{$$} && [name=A2]\MyBox{$u$}&& [name=A3]\MyBox{$$} \\
  &&[name=A4]\MyBox{$$}&& [name=A5]\MyBox{$$}   \\
  &[name=A6]\MyBox{$$} && [name=A7]\MyBox{$v$}&& [name=A8]\MyBox{$$}
  \psset{nodesep=0pt,arrows=-}
  \ncline{A1}{A2}
  \ncline{A2}{A3}
  \ncline{A1}{A4}
  \ncline{A1}{A6}
  \ncline{A3}{A8}
  \ncline{A3}{A5}
  \ncline{A4}{A5}
  \ncline{A4}{A6}
  \ncline{A5}{A8}
  \ncline{A6}{A7}
  \ncline{A7}{A8}
  \ncarc[arcangle=25,linecolor=red,linestyle=dashed]{A2}{A7}
 \end{psmatrix}
 &
 \begin{psmatrix}[colsep=0.5,rowsep=0.2,mnode=circle]
  &[name=A7]\MyBox{$v$} &&& [name=A1]\MyBox{$$} \\
  &&&[name=A5]\MyBox{$$} && [name=A2]\MyBox{$$} \\
  && && [name=A6]\MyBox{$u$}  \\
  &&&[name=A4]\MyBox{$$} && [name=A3]\MyBox{$$}
  \psset{nodesep=0pt,arrows=-}
  \ncline{A1}{A2}
  \ncline{A2}{A3}
  \ncline{A3}{A4}
  \ncline{A4}{A5}
  \ncline{A1}{A5}
  \ncline{A2}{A6}
  \ncline{A3}{A6}
  \ncline{A5}{A6}
  \ncline{A1}{A7}
  \ncline{A4}{A7}
  \ncarc[arcangle=45,linecolor=red,linestyle=dashed]{A6}{A7}
 \end{psmatrix} \\
 (a) & (b) & (c)
 \end{tabular}
  \caption{Implicit edges in 3-coloring of 3-chromatic graphs.}
  \label{impl1}
\end{centering}
\end{figure}

Also, in the graph (c) of Fig.~\ref{impl1} exist and implicit edge between
the vertices \{$u$,$v$\}, because its neighborhood is a 5-cycle and a
5-cycle is not 2-colorable. In this example we show a more generalizable implication
structure, showing that if two vertices conforms the neighborhood
of a subgraph H with $\chi(H)=3$ then; these two vertices can not
be colored with the same color because its adjacency is not
2-colorable.

\begin{theorem}Implicit Edges and Independent Sets.\\
\label{IE2} Given a k-chromatic graph G, the subset $\{i,j\}$ of vertices of G is an implicit-edge if and only
if there is no independent set $S$ of $G-ij$ containing the vertices $\{i,j\}$ such that the removal of $S$ from G results in a subgraph H of G with a chromatic number less than k.

\begin{proof}
Given a k-chromatic graph G and let $i$ and $j$ be to vertices of $G$:
\\
1) Assume that $\{$i$,$j$\}$ is an implicit-edge.

Let us suppose that by removing an independent set $S$ of $G-ij$ containing the vertices $i$ and $j$ the resulting graph $H=G-S$ is (k-1)-colorable. Then we can find a (k-1)-coloring of the graph $H$ and restore all the vertices in $S$ and assign them the color k and thus there will be a k-coloring of the graph $G-ij$ where vertices $i$ and $j$ are colored with the same color and this contradicts the assumption that $\{$i$,$j$\}$ is an implicit-edge.
\\
\\
2) Assume that $\{$i$,$j$\}$ is not an implicit-edge.

Then by definition there is a k-coloring of $G-ij$ such that the vertices $\{$i$,$j$\}$ are colored with the same color and then there will be an independent set $S$ of $G-ij$ containing the vertices $\{$i$,$j$\}$ (i.e. a color class) such that the removal of $S$ from G will result in a subgraph H of G with a chromatic number less than k. 
\end{proof}
\end{theorem}

As it is known another way of defining the graph coloring problem is the task of partitioning the graph into different independent sets and thus the chromatic number is then defined as the minimum number of partitions that can be achieved. Here the interesting fact is the relation of this definition and the Theorem~\ref{IE2} that yields the definition of a critical independent set.
\begin{definition}Critical Independent Set.\\
Given a k-chromatic graph G, we say that an independent set $S$ is
critical if and only if $\chi(G-S)=(k-1)$, that is, $S$ is a color
class of some k-coloring of G.
\end{definition}
From this definition follows that critical independent sets can not contain implicit-edges.

\begin{theorem}Implicit Edge invariant.\\
\label{ieinvariant}
If there are implicit-edges in a k-chromatic graph then by removing any critical independent set the resulting graph is (k-1)-chromatic with at least the same implicit-edges.
\begin{proof}
Given a k-chromatic graph G and let i and j be to vertices of G where $\{i,j\}$ is an implicit-edge:
\\
Remove a critical independent set of G, let us suppose that now $\{i,j\}$ is not an implicit-edge thus find a coloring where $\{i,j\}$ receive the same color, then restore the critical independent set with all vertices colored the same. The graph G is k-colored and $\{i,j\}$ has the same color and this contradicts the fact that $\{i,j\}$ is an implicit-edge of G.
\end{proof}
\end{theorem}

\begin{theorem} Implicit Edge (alternative definition on bipartite subgraphs).\\
\label{IE3}

Given a k-chromatic graph G. A subset $\{i,j\}$ of vertices of G
is an implicit-edge \textbf{if~and~only~if} $\{i,j\}$ is an
implicit-edge of $H$ for all subgraph $H<G$ such that:
$$H=G-\{S_1,\ldots,S_{k-2}\}~~;~~ \mbox{for any combination of
cardinality $(k-2)$ of $S_i$},$$ being each $S_i$ any critical
independent set of G (a color class of some k-coloring of G) and H
a bipartite graph, that is, vertices $\{i,j\}$ have distinct color
in every 2-coloring of H, that is, there is an odd length path
between $\{i,j\}$ after any arbitrary deletion of $(k-2)$ color
classes of G.

\begin{proof}$ $\\
Given a k-chromatic graph G if $\{i,j\}$ is an edge of G delete
this edge from G.
\\
\\
i) Let the subset $\{i,j\}$ of vertices of G be an implicit-edge.
Then by Theorem~\ref{iekch} in all k-colorings of G will exists a
kempe-chain $\Omega_{ij}$ of length $l\geq 3$ with colors
$(c(i),c(j))$ where $(c(i)\neq c(j))$ and vertices $\{i,j\}$ in
its extremes thus by deleting $(k-2)$ color classes in any
k-coloring of G you get a bipartite subgraph H of G (the
respective kempe-chain) in which vertices $\{i,j\}$ receive
different color, that is, $\{i,j\}$ is an implicit edge of H.
\\
\\
ii) Let the subset $\{i,j\}$ of vertices of G not be an
implicit-edge. Then simply find a k-coloring of G such that
vertices $\{i,j\}$ are colored with the same color, thus you can
delete $k-2$ color classes and obtain a 2-colorable subgraph H of
G in which vertices $\{i,j\}$ are colored with the same color and
hence $\{i,j\}$ is not an implicit-edge of some H.
\end{proof}
\end{theorem}

It is necessary to emphasize the importance of Theorem~\ref{IE3}.
This theorem implies that we can explain Implicit Edges in
k-chromatic graphs in terms of induced bipartite subgraphs. In
simple words: if $\{i,j\}$ is an implicit-edge of a k-chromatic
graph G then by deleting (k-2) arbitrary critical independent sets
(color classes of G) we obtain a bipartite graph $H<G$ where
$\{i,j\}$ is an implicit-edge too, and, if $\{i,j\}$ is not an
implicit-edge of G, then would be a combination of cardinality
(k-2) critical independent sets (color classes of G) such that
$H=G-\{S_1,\ldots,S_{k-2}\}$ is at most bipartite (2-colorable)
and $\{i,j\}$ is not an implicit-edge of $H$, this is due to
kempe-chains definition of implicit-edge.

\begin{definition}Explicit Neighborhood.\\
\label{TN}
Given a k-chromatic graph G, The Explicit Neighborhood
$EN_{k}(v)$ of one $v\in V$ is the complete set of vertices that
are joined by implicit or explicit edges to a vertex $v$ of $G$.
\end{definition}

This concept extends the notion of neighborhood to include the
implicit edges. As an example, we can see in Figure \ref{impl1}
(a) that the neighborhood of vertex \{2\} is: $\{1,4\}$ but its
Explicit Neighborhood for every 3-coloring of G: $EN^{3}(v=2)$ is
$\{1,4,5\}$.

\begin{definition}Explicit Graph.\\
Given an integer k and a graph G, an Explicit Graph is a graph $G'_k$ such that $G'_k$ does not have
implicit edges for a k-coloring problem of G.
\end{definition}

We can define the operator $G'=\Psi^{k}(G)$ as the operator which
represents a procedure that transform the graph $G$ into its
explicit representation, that is, $\Psi^{k}(G)$ copy the graph G
into $G'$ and add all the implicit edges of G to $G'$ for a
k-coloring problem.

It is important to note that critical graphs are explicit graphs,
that is, critical graphs does not have implicit-relations, this
fact will be shown in the section about critical graphs.

\subsection{Implicit Identities}
Implicit-Identities are in most cases the opposite-equivalent of Implicit-Edges, so there will be an almost one to one collection of theorems and lemmas which are equivalent. While implicit-edges are defined as a relation in which vertices receive always different color implicit-identities are the opposite relation, that is, the vertices receive always the same color.

\begin{definition}Implicit Identity. \\
\label{II} Given a k-chromatic graph G, we say that the subset
$\{i,j\}$ of vertices of G where $\{i,j\}$ is not an edge of G is
an implicit-identity if and only if in every k-coloring of G the
subset $\{i,j\}$ is colored with the same color.
\end{definition}

Note that contrary to implicit-edges, if we add and edge $e=\{i,j\}$ to a graph G the resulting graph $G+e$ will not be k-colorable and $e$ will be a critical edge of $G+e$.

\begin{theorem}Implicit Identities and Independent Sets.\\
\label{II2} Given a k-chromatic graph G, the subset $\{i,j\}$ of vertices of G where $\{i,j\}$ is not an edge of G is
an implicit-identity \textbf{if and only if there is no}
independent set $S$ in G containing the vertex $i$ or containing
vertex $j$ but not both such that, the complete removal of $S$
from G, results in a subgraph $H<G$ with chromatic number less
than k.
\begin{proof}
Given a k-chromatic graph G and let i and j be to vertices of G:
\\
1) Assume that $\{i,j\}$ is an implicit-identity.

Let us suppose that by removing an independent set $S$ of $G$ containing the vertex i or j (but not both) the resulting graph $H=G-S$ will be (k-1)-colorable. Then we can find a (k-1)-coloring of $H$ and restore all the vertices of $S$ assigning them the color k, thus there will be a k-coloring of $G$  where the vertices i and j will be colored with distinct color and this contradicts the assumption that $\{i,j\}$ is an implicit-identity.
\\
\\
2) Assume that $\{i,j\}$ is not an implicit-identity.

Then by definition there is a k-coloring of $G$ such that the independent set $\{i,j\}$ is colored with different color. Then there is an independent set $S$ of $G$ containing the vertex i or j (but not both) (i.e. a color class) such that the removal of $S$ from G will result in a subgraph H of G with a chromatic number less than k.
\end{proof}
\end{theorem}

\begin{theorem} Implicit Identity Invariant. \\
\label{ieinvariant}
If there is an implicit-identity $\{i,j\}$ in a k-chromatic graph then by removing any critical independent set not containing $\{i,j\}$ the resulting graph is (k-1)-chromatic with $\{i,j\}$ being an implicit-identity.
\begin{proof}
Given a k-chromatic graph G, let i and j be to vertices of G such that $\{i,j\}$ is an implicit-identity.
\\
Remove a critical independent set $S$ of G not containing i nor j and find a (k-1)-coloring of $G$ where $\{i,j\}$ receive different color. Then restore the critical independent set $S$ with all vertices colored the same. Thus the graph G is now k-colored and $\{i,j\}$ have different colors, but this contradicts the fact that $\{i,j\}$ is an implicit-identity of G.
\end{proof}
\end{theorem}

\begin{theorem}Implicit Identity (alternative definition on bipartite subgraphs).\\
\label{II3} Given a k-chromatic graph G(V,E), the pair
$\{i,j\}\in V(G)$ is an implicit-identity of G \textbf{if and only
if} $\{i,j\}\in V(H)$ is an implicit-identity of $H$ for all $H<G$
such that:
\begin{displaymath}
H=G-\{S_1,...,S_{k-2}\}~~;~~ \mbox{for any combination of
cardinality $(k-2)$ of $S_i$ },
\end{displaymath}
being each $S_i$ any critical independent set in G (i.e. a color
class) not including vertices $\{i,j\}$ and H a bipartite graph.
That is, vertices $\{i,j\}\in V(H)$ have the same color in every
proper coloring of H, that is, there is an even length path
between $\{i,j\}$ after any arbitrary deletion of (k-2) color
classes of G not containing vertices $\{i,j\}\in V(G)$.
\begin{proof}$ $\\
Given a k-chromatic graph G with a pair of vertices $\{i,j\}$ where $\{i,j\}$ is not an edge, then:
\\
\\
i) Let the subset $\{i,j\}$ of vertices of G be an implicit-identity.
Then by Theorem~\ref{iikch} in all k-colorings of G will exists at least
kempe-chain $\Omega_{ij}$ of length $l\geq 2$ with colors
$(c(i),c(j))$ where $(c(i)=c(j))$ and vertices $\{i,j\}$ in
its extremes thus by deleting $(k-2)$ color classes in any
k-coloring of G you get a bipartite subgraph H of G (the
respective kempe-chains) in which vertices $\{i,j\}$ receive
the same color, that is, $\{i,j\}$ is an implicit edge of H.
\\
\\
ii) Let the subset $\{i,j\}$ of vertices of G not to be an
implicit-identity. Then simply find a k-coloring of G such that
vertices $\{i,j\}$ are colored with different color, thus you can
delete $k-2$ color classes and obtain a 2-colorable subgraph H of
G in which vertices $\{i,j\}$ are colored with different color and
hence $\{i,j\}$ is not an implicit-identity of some H.
\end{proof}
\end{theorem}

It is necessary to emphasize the importance of
Theorem~\ref{II3}. This Definition implies that we can explain
Implicit Identities in k-chromatic graphs in terms of induced
bipartite subgraphs. In simple words: if $\{i,j\}$ is an
implicit-identity of a k-chromatic graph G then by deleting (k-2)
arbitrary critical independent sets (color classes of G) no
containing $\{i,j\}\in V(G)$ we obtain a bipartite graph $H<G$
where $\{i,j\}$ is an implicit-identity too, and, if $\{i,j\}$ is
not an implicit-identity of G, then would be a combination of
cardinality (k-2) critical independent sets (color classes of G)
such that: $$H=G-\{S_1,...,S_{k-2}\}~~~,$$ is at most bipartite
(2-colorable) and $\{i,j\}$ is not an implicit-identity of $H$.

This is important because implicit-identities in bipartite graphs
are defined very easy: every pair of vertices with an even length
path between them are implicit identities, i.e.
$i\bullet$---$\bullet$---$\bullet j$ then $\{i,j\}$ is an
implicit-identity.

\subsection{Chromatic Polynomials} A fundamental field of
research in graph coloring theory are the Chromatic Polynomials
$P(G,k)$, see~\cite{Dong2005}.
\\
\\
The chromatic polynomial $P(G,k)$  for a given graph is a
polynomial which encodes the number of different k-colorings of a
graph G, thus, we can denote the number of proper k-colorings of a
graph G as $P(G,k)$. The Chromatic Polynomials where first used by
Birkhoff and Lewis \cite{Birkhoff12,Birkhoff46} in their attack on
the four-color theorem.

One of the defining features of the chromatic polynomial $P(G,k)$
also called the Fundamental Reduction Theorem
(F.R.T)~\cite{Dong2005}(pp.4--6) of a graph G is that it satisfies
the deletion-contraction relation:

If $\{e\}$ is an edge of G joining vertices $\{x\}$ and $\{y\}$,
let $G-e$ be the graph G with $\{e\}$ deleted and let $G/e$ be the
graph resulting from the contraction of vertices $\{x\}$ and
$\{y\}$.
\\
Then:
\begin{equation}
\label{eq1} P(G,k)=P(G-e,k)-P(G/e,k)
\end{equation}

Also, the next equation is part of the so called F.R.T. Let
$\{x\}$ and $\{y\}$ be two non-adjacent vertices in a graph G and
let $G+xy$ be the graph with a new edge~$e=\{x,y\}$.
\\
Then:
\begin{equation}
\label{eq2} P(G,k)= P(G + e,k) +  P(G/e,k)
\end{equation}
A direct consequence of Eq.~\ref{eq1}, Eq.~\ref{eq2} and the
definition of implicit relations is presented as two new theorems~(\ref{ie:tdcr},\ref{ii:tdcr}):

\begin{theorem} Fundamental Reduction Theorem under implicit-edges.\\
\label{ie:tdcr} If G is a k-chromatic graph and $e= \{x,y\}$ is an
implicit edge of G (drawn~or~not).
\\
Then $P(G/e,k)=0$, that is, the graph  $( G/e ) $ is not
k-colorable.
\begin{proof}{of Theorem \ref{ie:tdcr}}\\
Being \{e\} an implicit-edge of a k-chromatic graph G: $P(G,k) =
P(G-e,k)$ or respectively (drawn or not) $P(G,k) = P(G+e,k)$ thus:
\begin{displaymath}
P(G,k) = P(G-e,k) \Longrightarrow  P(G/e,k)=0 \quad \mbox{due to
Eq.~\ref{eq1}}.
\end{displaymath}
\begin{displaymath}
P(G,k) = P(G+e,k) \Longrightarrow  P(G/e,k)=0 \quad \mbox{due to
Eq.~\ref{eq2}}.
\end{displaymath}
\end{proof}
\begin{corollary} of theorem~\ref{ie:tdcr}\\
Given a minor-closed class $\;\mathcal{G}$ of graphs (i.e. if $G
\in \mathcal{G}$ and $G \geq H$, then $H\in \mathcal{G}$).
\\
Being G a k-chromatic graph in $\mathcal{G}$ and \{e\} a drawn
implicit edge of G. If we contract the edge \{e\} the resulting
graph is a (k+1)-chromatic graph in $\mathcal{G}$. Also, in
general any contraction of an implicit edge (drawn or not) will
result in a (k+1)-chromatic graph.
\end{corollary}
\end{theorem}

And in the same sense for the implicit-identities case:

\begin{theorem} Fundamental Reduction Theorem under implicit-identities.\\
\label{ii:tdcr} If G is a k-chromatic graph and $e= \{x,y\}$ is an
implicit-identity of G.
\\
Then $P(G,k) = P(G/e,k)$.
\begin{proof}{of Theorem \ref{ii:tdcr}}\\
Being \{e\} an implicit-identity of a k-chromatic graph G: $P(G+e,k)=0$ thus:
\begin{displaymath}
P(G+e,k)=0 \Longrightarrow  P(G,k) = P(G/e,k) \quad \mbox{due to Eq.~\ref{eq2}}.
\end{displaymath}

\end{proof}
\end{theorem}

\subsection{Implicit Relations in planar graphs}
An important type of graphs are the planar graphs, a graph is
called planar: if it can be drawn in a plane without edges
crossings. It is important to show that 3-coloring planar graphs
remains NP-Complete \cite{GJ79}. Planar graphs poses the next
important properties:
\begin{description}
    \item[a)] Every planar graph is four colorable
              \cite{ah1,ah2,ah3,Birkhoff,Robertson,RSST}.
    \item[b)] Every planar graph without triangles is
              3-colorable\cite{JT95,gb63,gro59}.
    \item[c)] Planar graphs are closed under edge contraction, that is, every
              edge contraction of a planar graph results in a new planar graph.
\end{description}

The most important property for our purposes is the property (a)
due to the celebrated four-color-conjecture that now is called
after \cite{ah1,ah2,ah3} the four color theorem. Given the fact
that every planar graph is four colorable there are some
interesting questions:
\begin{quote}
- There are implicit edges in 4-coloring of planar graphs?
\\
- Are the proofs of \cite{ah1,ah2,ah3,Birkhoff,Robertson,RSST}
affected in some way by the existence of implicit edges?.
\end{quote}

\begin{figure}[tbh]
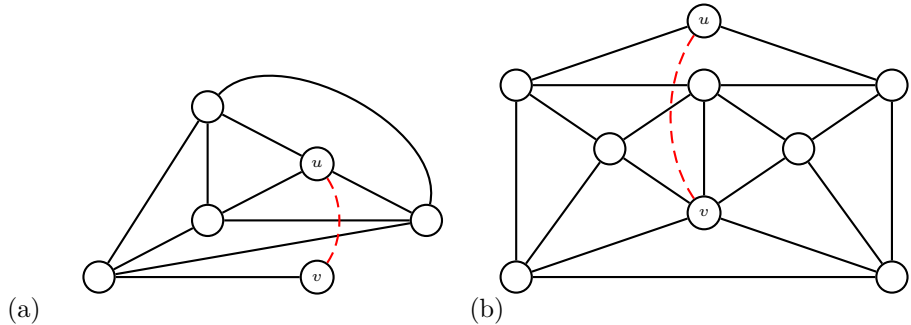

\begin{centering}
  \begin{tabular}{ll}
   \begin{psmatrix}[colsep=1.0,rowsep=0.3,mnode=circle]
  &&[name=A2]\MyBox{$$} \\
  &&&[name=A1]\MyBox{$u$}\\
  &&[name=A3]\MyBox{$$} &&[name=A4]\MyBox{$$}\\
  &[name=A5]\MyBox{$$} &&[name=A6]\MyBox{$v$}

  \psset{nodesep=0pt,arrows=-}
  \ncline{A1}{A2}
  \ncline{A2}{A3}
  \ncline{A3}{A4}
  \ncline{A4}{A1}
  \ncline{A1}{A3}
  \ncarc[arcangle=75,linecolor=black,linestyle=solid]{A2}{A4}
  \ncline{A5}{A2}
  \ncline{A5}{A3}
  \ncline{A5}{A4}
  \ncline{A5}{A6}
  \ncarc[arcangle=35,linecolor=red,linestyle=dashed]{A1}{A6}

 \end{psmatrix}
 &
 \begin{psmatrix}[colsep=0.4,rowsep=0.4,mnode=circle]
  &               &&&& [name=A0]\MyBox{$u$} \\
  &[name=A1]\MyBox{$$} &&&& [name=A2]\MyBox{$$} &&&& [name=A3]\MyBox{$$} \\
  &&&[name=A4]\MyBox{$$}&&&& [name=A5]\MyBox{$$}   \\
  & &&&& [name=A7]\MyBox{$v$} \\
  &[name=A6]\MyBox{$$} &&&& &&&& [name=A8]\MyBox{$$}
  \psset{nodesep=0pt,arrows=-}
  \ncline{A0}{A1}
  \ncline{A0}{A3}

  \ncline{A1}{A2}
  \ncline{A2}{A3}

  \ncline{A1}{A4}
  \ncline{A4}{A6}

  \ncline{A3}{A5}
  \ncline{A5}{A8}

  \ncline{A2}{A4}
  \ncline{A2}{A5}
  \ncline{A2}{A7}
  \ncline{A4}{A7}
  \ncline{A5}{A7}

  \ncline{A7}{A6}
  \ncline{A7}{A8}
  \ncline{A6}{A9}

  \ncline{A1}{A6}
  \ncline{A3}{A8}
  \ncline{A6}{A8}

  \ncarc[arcangle=-35,linecolor=red,linestyle=dashed]{A0}{A7}
 \end{psmatrix} \\
 (a) & (b)
 \end{tabular}
  \caption{Implicit edges in 4-chromatic graphs.}
  \label{planar4c}
\end{centering}
\end{figure}

The first question is affirmative, there are implicit edges in
4-coloring of planar graphs as can be seen in Figure
\ref{planar4c}.

The main observation from Figure \ref{planar4c} is that the
implicit edge can not be drawn without edge line crossings, it is
straightforward to proof the fact that if implicit edges in
4-coloring of planar graphs could be drawn without edge line
crossings, then the theorems in~\cite{ah1,ah2,ah3,Birkhoff,Robertson,RSST} are false due to the
corollary of Theorem~\ref{ie:tdcr}.

\subsection{Implicit Relations and Kempe chains}
\label{sec:14}

Kempe chains were first used by Alfred Kempe in his attempt of a
proof of the four colors theorem~\cite{kempe}. Although his proof
was not valid his method of color chains remained as a useful
tool, indeed it is used in modern proofs like new computer based
proofs of the four colors theorem see \cite{ah3,Robertson}. Next
we provide a simple definition of Kempe Chains.

Let G be a graph with vertex set V, and we provide a proper
coloring function:
\begin{displaymath}
\varphi : V(G)\rightarrow C~~~,
\end{displaymath}
where C is a finite set of colors, containing at least two
distinct colors $c_1$ and $c_2$. If $\{v\}$ is a vertex with color
$c_1$, then the $(c_1,c_2)$-Kempe chain of G is the maximal
connected subset of V which contains $\{v\}$ and whose vertices
are all colored either $c_1$ or $c_2$, i.e. a simple bipartite
component of G colored with colors $c_1$ and $c_2$.

One of the useful characteristics of Kempe chains is the well
known process of flipping. This process consist in taking a kempe
chain $\Omega$ of colors $c_1$ and $c_2$ and change all vertices
of $\Omega$ with color $c_1$ to $c_2$ and all vertices with color
$c_2$ to $c_1$, indeed this was the basic mechanism of the failed
proof the four color theorem by Kempe.

Kempe chains could be easily related to both implicit-edges and
implicit-identities but the most strong relation is for the
implicit-identities case. The relation between kempe-chains an implicit-relations is important since the existence of bipartite chains and path is a useful tool to define some structural properties of some special classes of graphs.

\begin{theorem} Kempe chains and implicit-identities.\\
\label{iikch}
Let G be a k-chromatic graph with an implicit-identity between a
pair of vertices $\{x,y\}$ of G, then in every proper k-coloring
of G will exist at least $(k-1)$ kempe chains with colors
$(c(x,y),c_i)$ with the vertices $\{x,y\}$ in its extremes.

\begin{proof} $ $ \\
Otherwise we can simply flip the color of vertex $\{x\}$ with some
color of an adjacent vertex and will exists a proper coloring of G
where $\{x,y\}$ are colored with distinct color.
\end{proof}
\end{theorem}

\begin{theorem} Kempe chains and implicit-edges.\\
\label{iekch}
Let G be a k-chromatic graph with an implicit-edge
between a pair of vertices $\{x,y\}$ of G, then in every proper
k-coloring of G will exist at least a Kempe chain with colors
$(c(x),c(y))$ where $(c(x)\neq c(y))$ and vertices $\{x,y\}$ in
its extremes.

\begin{proof} $ $ \\
Otherwise (by deleting a possible edge $\{x,y\}$) we can simply
flip the color of vertex $\{x\}$ with some color of an adjacent
vertex with the color of vertex $\{x\}$ and will exists a proper
coloring of G where $\{x,y\}$ are colored with the same color.
\end{proof}
\end{theorem}

We should recall that implicit-edges could belongs to the set of
edges or not, thus in the case of Lemma~\ref{iekch}, if the
implicit-edge belongs to the set of edges the kempe-chain will
have an odd cycle.

\subsection{Implicit Relations in Critical Graphs}

Critical graphs where first studied by Dirac~\cite{dirac52,dirac52c,dirac53}. An element $\{x\}$ of a
graph G is called critical if $\chi(G-x)<\chi(G)$. If every element $\{x\}$ (vertex or edge) of a graph G is critical and
$\chi(G)=k$ we say that G is k-critical, examples of critical graphs are the complete graphs $K_k$. If G is critical and
different from $K_k$ then G has $n\geq k+2$ vertices.

A k-chromatic graph G is called double-critical if for every edge $\{x,y\}$ of G the graph $G-x-y$ is $(k-2)$-colorable, and example
of a double-critical graph is the complete graph $K_k$.

It is easy to see that if a k-chromatic graph G is double-critical then it is k-critical, otherwise by deleting any non critical
vertex $\{x\}$ of G one can delete a neighbor $\{y\}$ of $\{x\}$ and then $G-x-y$ is not $(k-2)$-colorable.

In every k-critical graph G different from $K_k$ there is at least a triplet $t= \{x,y,z\}$ of vertices of $G$ such that $\{x,y\}$ is
an edge of $G$, $\{y,z\}$ is an edge of $G$ but $\{x,z\}$ is not an edge of $G$ (i.e. a simple path $x\bullet$---$\bullet
y$---$\bullet z$), otherwise the graph is a complete graph $K_k$.

Also k-critical graphs posses the next well known properties:
\begin{itemize}
    \item G has only one component.
    \item G is finite
    \item Every vertex is adjacent to at least k - 1 others.
    \item If G is (k-1)-regular, meaning every vertex is adjacent to exactly
          k - 1 others, then G is either $K_k$ or an odd cycle.
    \item $2|E| \geq (k - 1)|V| + k - 3$, indeed, $2|E| \geq (k - 1)|V| + [\frac{k - 3}{k2 - 3}]|V|$
    \item $|V|\neq k+1$
\end{itemize}
\noindent being k the chromatic number.
\\
\\
One interesting characteristic of implicit-relations in critical graphs is that there are no implicit relations in critical graphs. This could be seen in the next theorem:

\begin{theorem}\label{noircritical} there are no implicit relations in critical graphs
\begin{proof}$ $ \\
\noindent i) Let us suppose a k-critical graph G such that G has an implicit-edge $e=\{a,b\}$. If $\{e\} \in E(G)$ delete edge
$\{e\}$ from G. Since G-ab is still k-chromatic there is a contradiction with the assumption that G is k-critical. If $\{e\} \notin E(G)$ the graph with vertices a and b removed ($G-a-b$) remains k-chromatic due to $e=\{a,b\}$ is an implicit-edge which contradicts the assumption that G is k-critical
\\
\\
ii) Let us suppose a k-critical graph G such that G has an implicit-identity $\{a,b\}$. Since $\{a,b\}$ is an implicit-identity the graph $G-a$ remains k-chromatic which is absurd due to G is critical.
\end{proof}
\end{theorem}
Indeed, now we can say also that in a critical graph all vertices, edges and independent sets are critical.

Although in critical graphs there are no implicit relations there are deep connections between critical graphs, critical elements an implicit relations.

Let us consider a k-critical graph G. It is easy to proof by means of implicit-relations that any vertex-identification in a
critical graph results in a k-chromatic graph. Also maybe the most interesting connection is in the study of the k-critical graphs
$G^{-}$ which are k-critical graphs with one edge removed. The study of the critical graphs $G^{-}$ is in deep connection with many interesting problems such as the four-colors theorem, the Hadwiger Conjecture, the double-critical graphs problem and so on.

It is well known that if some edge $e=\{a,b\}$ of a k-chromatic graph G is critical then the removal of this edge results in a
(k-1)-chromatic graph $G^{-}$ where all proper (k-1)-colorings of $G^{-}$ assigns the same color to vertices $\{a,b\}$, thus
vertices $\{a,b\}$ in $G^{-}$ are an implicit-identity and then will be at least (k-2)-kempe-chains between $\{a,b\}$ in any
proper (k-1)-coloring of $G^{-}$.

Even more, another interesting result is that any graph with at least one critical vertex $(z)$ should have for every implicit-edge $(u,v)$ at least one of its vertices $u$ or $v$ adjacent to $z$. In the same way any graph with at least one critical vertex $(z)$ should have for every implicit-identity $(u,v)$ all of its vertices $u$ and $v$ adjacent to $z$.

\section{Complexity}

\subsection{Polynomial 3-coloring}
In the next theorem it is shown that if there is any polynomial
algorithm $\wp$ for determining implicit-edges then there is a
polynomial algorithm to determine 3-colorability and then $P=NP$.

\begin{theorem}polynomial 3-coloring (I).\\
\label{t4}

Given a Graph $G(V,E)$, if there is a polynomial algorithm that
determine if two vertices $\{u,v\}$ are an implicit edge for a
3-coloring problem, then 3-coloring is in P.
\begin{proof}{of theorem \ref{t4}}
\\
Given a graph G and a polynomial algorithm $\wp_{ie}$ that
determine if two vertices $\{u,v\}$ are an implicit edge.
\\
\\
For every pair of vertices $\{u,v\}$ of G try $\wp_{ie}(u,v)$:
\\
$- \quad \quad if$ $\wp_{ie}(u,v)$, add this edge to G.
\\
\\
1.- if the resulting graph is a complete graph then G is not
3-colorable.
\\
2.- if the resulting graph is not a complete graph then G is
3-colorable.
\end{proof}
\end{theorem}

Also we can use together the definitions of implicit-edges and
implicit-identities in a more simpler algorithm.

Let $\wp_{ie}$ be a polynomial algorithm that determine if two
vertices $\{u,v\}$ are an implicit edge.

Let $\wp_{ii}$ be a polynomial algorithm that determine if two
vertices $\{u,v\}$ are an implicit identity.

For some pair of vertices $\{u,v\}$ of G try: $$[\wp_{ie}(u,v)
\land \wp_{ii}(u,v)]~~~,$$

if $\left[\wp_{ie}(u,v) \land \wp_{ii}(u,v)\right]=true$ then
there is a contradiction and thus the graph is not 3-colorable,
otherwise the graph is 3-colorable.

\subsection{NP and CoNP}

We have shown that the decision problem of implicit-edges and
implicit identities in general belongs to the class CoNP. So, Let
G be a graph G and let the pair $\{a,b\}$  of G be not joined by
an edge in G, the complete relations to complexity classes can be
seen in Table~\ref{tablecomplexity}.

\begin{table}[tbh]
\label{tablecomplexity}
\begin{centering}
\begin{tabular}{|l|r|r|}
  \hline
  DECISION PROBLEM & YES & NO  \\
  \hline
  Implicit-edge         & CoNP & NP \\
  Implicit-identity     & CoNP & NP \\
  no implicit relation  & NP & CoNP \\
  k-colorability of G   & NP & CoNP \\
  \hline
\end{tabular}
\caption{Complexity table of colorability}
\end{centering}
\end{table}

\section{Conclusions}
We have presented a theory about the implication structure in graph coloring. Possible applications of such theory can be
expected in the discovery of hidden relations in a wide range of scientific disciplines. The main contribution presented in this
work is a model of hidden relations materialized as implicit-edges and implicit-identities, these relations can be interpreted in
physical and chemical models as hidden forces, hidden interactions, hidden reactions and hidden variables.

\end{document}